\documentstyle[12pt,twoside,numinsec]{siamltex}
\input psfig
\input{amssym.def}
\input{amssym.tex}
\newtheorem{plmthm}{Theorem \ref{PLM-follows-ray}} 
\newtheorem{defn}[theorem]{Definition}
\newenvironment{algorithm}[1]%
    {\pagebreak[3]

     \newcommand{\stmt}{\nopagebreak[2]\item}
     \par\vskip 5pt plus 10pt
     \centerline{\tt#1}
     \medskip
     \nopagebreak[4]
     \begin{itemize}
            \sf\small%
            \baselineskip=1.1\baselineskip %
            }%
    {\end{itemize}
     \par\pagebreak[3]\vskip 5pt plus 10pt}

\newcommand{\thmref}[1]{Theorem~\ref{#1}}

\newcommand{\lemref}[1]{Lemma~\ref{#1}}
\newcommand{\defref}[1]{Definition~\ref{#1}}
\newcommand{\secref}[1]{Section~\ref{#1}}
\newcommand{\figref}[1]{Figure~\ref{#1}}

\newcommand{\dh}{\discretionary{}{}{}} 

\newcommand{\C}{{\Bbb C}}

\newcommand{\D}{{\Bbb D}}

\newcommand{\Order}[1]{\displaystyle{{\cal O}\negthinspace\({#1}\)}}
\newcommand{\Pd}{{\cal P}_d}
\newcommand{\Wedge}[1]{{\cal W}_{#1}}
\newcommand{\maps}{\rightarrow}

\newcommand{\union}{\cup}
\newcommand{\Union}{\bigcup}
\newcommand{\e}{{\rm e}}
\newcommand{\st}{\,\bigm|\,}
\renewcommand{\i}{{\rm i}}
\renewcommand{\(}{\left (}
\renewcommand{\)}{\right )}
\newcommand{\Label}[1]{%
	\label{#1}%
        }

\newcommand{\Cite}[1]{%
	\cite{#1}%
        }

\title	{Polynomial Root-Finding Algorithms\\
         and Branched Covers%
         \footnotemark[1]\ \footnotemark[2]}

\author {Myong-Hi Kim\footnotemark[3]\ \footnotemark[4]
         and Scott Sutherland\footnotemark[3]}

\begin{document}
\thispagestyle{empty}
\maketitle

\renewcommand{\thefootnote}{\fnsymbol{footnote}}

\footnotetext[1]{Received by the editors July 12, 1991.}
\footnotetext[2]{%
   An earlier version of this paper was circulated with the title 
   ``Parallel Families of Polynomial Root-Finding Algorithms''.}
\footnotetext[3]{%
   Institute For Mathematical Sciences,
   SUNY Stony Brook,
   Stony Brook, NY 11794-3660.
   email: {\tt myonghi@math.sunysb.edu}\ and
          {\tt scott@math.sunysb.edu}.}
\footnotetext[4]{%
   Part of this work was done while M.~Kim was at Bellcore,
   Morristown, NJ.}

\renewcommand{\thefootnote}{\arabic{footnote}}
\begin{abstract}
We construct a family of root-finding algorithms which combine 
knowledge of the branched covering structure of a polynomial
with a path-lifting algorithm for finding individual roots.
In particular, the family includes an algorithm
that computes an $\epsilon$-factorization of a polynomial of degree $d$
which has an arithmetic complexity of 
$\Order{d(\log d)^2|\log\epsilon| + d^2(\log d)^2}$.  
At the present time, this complexity is the best known in terms of the
degree.  
\end{abstract}



\section*{Introduction}
The problem of devising optimal methods for numerically approximating the
roots of a polynomial has been of interest for several centuries, and is far
from solved. 
There are numerous recent works on root-finding algorithms and their cost,
for example, the work of Jenkins and Traub \cite{JenkinsTraub:3StageShift},
Renegar \cite{Renegar}, Sch\"onhage \cite{Shonhage}, and Shub and Smale
\cite{ShubSmale:Complexity1,ShubSmale:Complexity2,Smale:bulletin85}. This
list is far from complete; the reader should refer to the aforementioned
papers as well as \Cite{DejonHenrici} for more detailed background.
The work in this paper is most closely related to that of Smale.

Our algorithm computes an approximate factorization of a given polynomial
(that is, it approximates all the roots).  In constructing it, we combine
global topological information about polynomials (namely, that they act as
branched covers of the Riemann Sphere) with a path-lifting method for finding
individual roots.  Utilizing this global information enables us to use fewer 
operations than applying the path-lifting method to each root sequentially.

\smallskip

Renegar's algorithm in \cite{Renegar} approximates all $d$ roots of a given
polynomial using $\Order{d^3\log d + d^2 (\log d) (\log|\log \epsilon|)}$
arithmetic operations in the worst case.  He has shown that the factor of
$\log|\log\epsilon|$ in the complexity is the best possible if one restricts
to the operations $+$, $-$, $\times$, and $\div$.  This algorithm has a
component (the Shur-Cohn algorithm) which requires exact computation and so
is not suitable for an analysis of bit complexity, that is, one which
accounts for rounding errors introduced by finite precision.  In
\cite{Pan:PolyZeros}, Pan gives a different algorithm which slightly
improves the complexity to $\Order{d^2\log d\log|\log\epsilon|}$; this
algorithm also operates effectively as a parallel algorithm.

Sch\"onhage \cite{Shonhage} gives an algorithm which produces an
$\epsilon$-factorization with a
bit complexity of  
 $\Order{d^3\log d +
  d^2|\log\epsilon|}\log(d|\log\epsilon|)\log\log(d|\log\epsilon|)$, 
via the ``splitting circle method''.  Note that the customary parameter for
bit length of the coefficients does not appear in the complexity.  This is
because, as Sch\"onhage states, for fixed degree $d$ and output precision
$\epsilon$, there is a number $s_0$ for which ``the input [coefficients]
$a_\nu$ can be restricted to complex integer multiples of $2^{-s_0}$ without
loss of generality.''
In \cite{Renegar}, it is stated that Sch\"onhage believes that, if exact
arithmetic is used, this method ``should yield a complexity bound [in
$\epsilon$] of $\Order{d^\alpha \log|\log\epsilon|}$, most likely with
$\alpha \le 3$.''

Smale's path lifting algorithm in \cite{Smale:bulletin85} approximates a single
root of the polynomial with a worst case arithmetic complexity of
$\Order{d(\log d)|\log\epsilon|+ d^2(\log d)^2 }$, and an average 
complexity of $\Order{d^2 +d|\log\epsilon|}$.  One good feature of this line
of work is that it is stable under round-off error.  In
\cite{Kim:BitComplexity}, Kim shows that if $f$ and $f'$ are computed with
relative error $10^{-3}$ until an approximate zero (see \secref{AZ-section})
is reached, then the algorithm behaves exactly the same.
A recent series of papers by Shub and Smale
\cite{ShubSmale:Bezout1,ShubSmale:Bezout2,ShubSmale:Bezout3,ShubSmale:Bezout4}
generalizes the path lifting algorithm to systems of homogeneous polynomials
in several variables. 
\smallskip

The algorithm presented here exploits the branched covering structure of a
polynomial to choose good starting points for a variant of Smale's
algorithm, and we obtain a worst case arithmetic complexity of  
$\Order{d(\log d)|\log\epsilon| + d^2(\log d)^2}$ to compute an
$\epsilon$-fact\-or\-ization. 
In a subsequent paper, we shall compute the bit complexity of this
algorithm.  Because of the stability mentioned in the previous paragraph and
the ability to exploit bounds on the variation of $f$ and $f'$, we hope to
achieve results comparable to Sch\"onhage's.  

At first glance, it may appear that our complexity results are inferior to
some of those above in terms of $\epsilon$.  However, in practice there is
usually a relationship between the degree $d$ and the desired precision
$\epsilon$; if we have $\epsilon \ge 2^{-d}$, then the complexity
of our algorithm compares favorably with all of those mentioned above.
Furthermore, our algorithm is quite simple to implement and is numerically
very stable. 

Our algorithm is suitable for some amount of parallelization, but has a
sequential component of $\Order{d + |\log\epsilon|}$ operations.  However, we
think of this algorithm as acting on $d$ points simultaneously, and
techniques which evaluate a polynomial at $d$ points (see
\Cite{BorodinMunro}, for example) are used to cut the cost involved.
Of course, the algorithm can be implemented on a sequential machine while
still taking advantage of these techniques.  In fact, evaluation of the
polynomial is the only point at which we at which we need to use asymptotic
estimates to achieve the stated complexity; the other places where we use
asymptotic estimates are only for ease of expositon.

The reader should also see the papers 
\Cite{BFKT:ParallelRealRoots,BenOrTiwari:SimpleAlgsForRealRoots,Neff:AlgsNCRootIsolation}
for fully parallel algorithms for solving polynomials with integer
coefficients.  In \cite{BFKT:ParallelRealRoots}, it is shown that if all
roots of the polynomial are real, this problem is in NC.  Neff extends this
result to allow complex roots in \cite{Neff:AlgsNCRootIsolation}.

\smallskip

This paper is structured as follows: In Chapter~\ref{prelims}, after some
background material, we recall the ``path lifting method'' of Smale and
present a version of the relevant theorem (our \thmref{PLM-follows-ray})
which improves the constants involved somewhat.  We then discuss how we can
exploit the branched-covering structure of a polynomial to choose initial
points for the algorithm, many of which will converge to roots. We close the
chapter with a brief explanation of how to construct families of algorithms
which locate $d/n$ roots at a time, for various values of $n$.

Chapter~\ref{algorithm} presents an explicit algorithm for a specific
family, which locates $d/2$ points at a time.  Our main theorem,
\thmref{main-theorem},  states that this algorithm always terminates with an
$\epsilon$-factorization of the input polynomial, and gives a 
bound on the number of arithmetic operations required in the worst case.
As a corollary, the algorithm can be used to locate all $d$ roots of the
polynomial to within $\epsilon$ with a complexity of 
$\Order{d^2(\log d)|\log\epsilon| + d^2(\log d)^2}$.
In the subsequent sections, each component of the algorithm is analyzed, and
the relevant lemmas are proven.  Finally, we tie all the components together
and prove the main theorem.

\bigskip
\section{Preliminaries} \label{prelims}

\subsection{Root and coefficient bounds} \Label{root-bounds}
Given a polynomial $\phi(z) = \sum_{i=0}^d a_i z^i$, with 
$a_i \in \C$, it is our goal to determine an approximate factorization of
$\phi$, that is, approximations $\hat\xi_i$ to the actual roots $\xi_i$ of
$\phi$ so that $\|\phi(z) - \prod(z - \hat\xi_i)\|<\epsilon$.  The norm we
shall use here is the max-norm,
that is, $\|\phi\| = \max |a_j|$. 
A related problem is to ensure that $|\xi_i - \hat\xi_i| < \epsilon'$; there
are well-known estimates giving the relationship between $\epsilon$ and
$\epsilon'$, so solving one problem essentially solves the other.

\medskip
In order to have an estimate on the complexity of a root-finding algorithm,
we need a compactness condition on the space of polynomials.  This can be
done either by placing conditions of the location of the roots or on the
coefficients; such bounds are interrelated.

Since our goal is to minimize a functional norm, it seems most natural to
place our conditions on the coefficients. We shall assume our input 
polynomial $\phi$ is an element of the family
$$\Pd(1) = \{ z^d + \sum_{j=0}^{d-1} a_j z^j,\ \hbox{with}\ |a_j| \le 1\}.$$
This is the same space as used by Smale and others
\Cite{ShubSmale:Complexity1,Smale:bulletin81,Smale:bulletin85}.
One can always transform an arbitrary polynomial into an element of
$\Pd(1)$: 
if $p(z) = \sum_{i=0}^d b_i z^i$, and $B=\max |b_j/b_d|^{1/(d-j)}$, then
$p(Bz)/B^d \in \Pd(1)$. 

One should not confuse this family with the degree $d$ polynomials whose
roots are in the unit disk, although unfortunately this space is also
often denoted by $\Pd(1)$ (for example, in \Cite{Friedman:ApproxZeros} and
\Cite{Renegar}).  

\medskip
There are a number of estimates which relate the coefficients of a
polynomial to a bound on the modulus of the zeros (see
\Cite{Henrici:AppliedComplexAnalysis} or \Cite{Marden:GeomPolys}, for
example).   Such an estimate is important to us, since although membership
in $\Pd(1)$ is not preserved under deflation (division of factors), bounds on
the modulus of the roots are.  We state one such bound here (Corollary~6.4k
of \Cite{Henrici:AppliedComplexAnalysis}):

\begin{lemma}
All the zeros of the polynomial $z^d + \sum_{j=0}^{d-1} a_j z^j$ lie within
the open disk with center 0 and radius
	$$2 \max_{0 \le j < d} | a_j |^{1/(d-j)}.$$
\end{lemma}

\noindent
As an immediate consequence, we see that the roots of a polynomial in
$\Pd(1)$ lie within $\D_2$.

\medskip
\subsection{Approximate zeros} \Label{AZ-section}
Our algorithm uses a path lifting method (see below) to get close to the
roots of our polynomial, and then uses the standard Newton's method to
further refine these approximations.  This is done because Newton's method
converges very quickly in a neighborhood of a simple root, but can fail 
for some initial points outside this neighborhood.
One of the authors \Cite{Sutherland:thesis} has shown how one
can guarantee convergence of Newton's method, but a bound on
the arithmetic complexity has not been computed.  Instead, we use the more
certain path lifting method as described in \secref{PLM-section}; this
allows an explicit computation of the complexity.

Following Smale \Cite{Smale:bulletin81}, we call a point $z_0$ an {\em
approximate zero} if Newton's method converges rapidly (that is,
quadratically) when started from $z_0$.  Such terminology is reasonable,
because given such a point, we can quickly obtain an approximation of a root
to arbitrary precision.

\begin{defn}      \Label{approx-zero}
  Let $f$ be a polynomial and let $z_n$ be the $n^{th}$ iterate under
  Newton's method of the point $z_0$, that is, 
  $z_{n} = z_{n-1} - f(z_{n-1})/f'(z_{n-1})$.
  Then we say that  $z_0$ is an approximate zero of $f$ if, for all $n>0$ we
  have  
	$$|z_n - \zeta| \le 8\(\frac{1}{2}\)^{2^n}|z_0-\zeta|,$$
  for some root $\zeta$ of $f$.
\end{defn}

Notice that this definition is never satisfied in the neighborhood of a
multiple root of $f$, since the convergence of Newton's method is
asymptotically linear there.  In our algorithm, we perturb the polynomial
slightly to ensure that we always have simple zeros.  Refer to
\secref{alg-statement} and \secref{PLM-iteration} for more details.

Kim \Cite{Kim:ApproxZeros} and Smale \Cite{Smale:NewtonEstimates} have
developed readily tested criteria for determining, based on the values of
the derivatives $f^{(k)}(z)$, when a point $z$ is an approximate zero. 
These can be extended to a much more general setting, namely for $f$ a
mapping between Banach spaces. The following is 
essentially Theorem~A of  \Cite{Smale:NewtonEstimates}:

\begin{lemma}   \Label{alpha}
  Let 
  	$$\alpha_f(z) =  \max_{k>1} \left|\frac{f(z)}{f'(z)}\right|
          \left|\frac{f^{(k)}(z)}{k!\,f'(z)}\right| ^{1/(k-1)}.$$
  If $\alpha_f(z) < \frac{1}{8}$, then $z$ is an approximate zero of $f$.
\end{lemma}

We will find the following also very useful.

\begin{lemma}   \Label{approx-zeros-if-radius-small}
  Let $f(z)$ be a polynomial and $z$ be a complex number so that
  $f'(z) \ne 0$, and let $R_f(z)$ be the radius of convergence of the branch 
  of the inverse $f^{-1}_{z}$ which takes $f(z)$ to $z$.  If
  	$$\frac{|f(z)|}{R_f(z)} < \frac{1}{10},$$
  then $z$ is an approximate zero of $f$.  Furthermore, if we have
  	${|f(z)|}/{R_f(z)} < {1}/{32},$
  then $\alpha_f(z) < 1/8$.
\end{lemma}

\begin{remark}
  If Smale's mean value conjecture holds (see
  \Cite{Smale:bulletin81,Tischler:CritPointsAndValues}), then the hypotheses
  of the lemma imply that $\alpha_f(z) < 1/8$.
\end{remark}

\begin{proof}
  The first result is a consequence of the proof of Theorem~4.4 of
  \Cite{Kim:ApproxZeros}, and the second is an immediate consequence of
  Corollary 4.3 of the same paper. This, in turn, uses the Extended Loewner's
  Theorem in \Cite{Smale:bulletin81}.
\qquad\end{proof}

\medskip
\subsection{The path lifting method}       \Label{PLM-section}
Here we review the path lifting method, which forms the core of our iteration
scheme.  This method is sometimes referred to as a ``generalized Euler
method'' or ``modified Newton's method''; we prefer the term ``path lifting
method'' as it is the most descriptive (to us, anyway).  This method appears
in the work of Steven Smale \Cite{Smale:bulletin85}, although the version we
present here is slightly different and we present another proof of the 
relevant theorem, which is quite simple.  It should be emphasized 
that the path lifting method, like Newton iteration, is an algorithm for
finding a {\em single} root of a polynomial; we discuss how to combine these
to find all roots in \secref{root-finding-families} below.

We think of a polynomial $f$ as a map from
the {\em source space} to the {\em target space}; that is,
$f:\C_{source}\maps\C_{target}$.  Given an initial value $z_0$ in the source
space, we connect its image $w_0 = f(z_0)$ to $0$ in the target space, and
then lift this ray under the proper branch of $f^{-1}$ to a path connecting
$z_0$ with a root $\zeta$ of $f$.  Of course, we don't explicitly know this
inverse, but if the path in the target space stays well
away from the critical values of $f$, the local inverse map $f^{-1}_{z_0}$ is
well-defined on a neighborhood of the ray.  Even if the path {\em does}
contain critical values, a local inverse can still be defined for some
$z_0$.  See \secref{branched-covers-sec}.

\begin{figure}[thb]
 \centerline{\psfig{figure=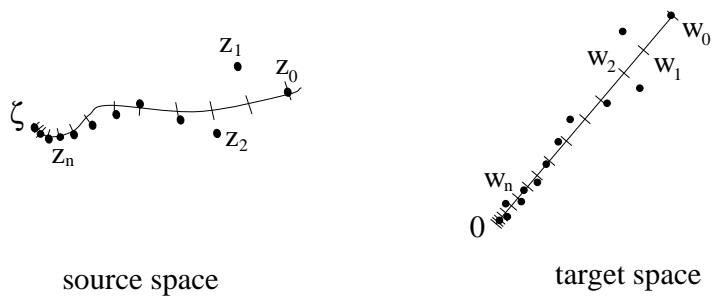,height=1.5in}}
 \caption{\Label{PLM-fig} The source and target spaces in the path lifting
         method.  In the source space, each $z_i$ is indicated by a black
         dot, and $f(z_i)$ is indicated by a black dot in the target space.
         Similarly, the $w_i$ are indicated by tick marks in the target
         space, and $f^{-1}(w_i)$ by ticks in the source space.} 
\end{figure}

The basic idea of the path lifting method is to take a sequence of
points $w_n$ along the ray in the target space, with $w_0 = f(z_0)$.  We then
construct a sequence of points $z_n$ in the source space so that
$f(z_n)$ is near $w_n$ in the target.  This is done using a single step
of Newton's method to solve $f(z)=w_n$ with $z_{n-1}$ as the starting point.
That is,
	$$z_n = z_{n-1} - \frac{f(z_{n-1})-w_{n-1}}{f'(z_{n-1})}.$$

This construction will converge to a root $\zeta$ in the source space if
there is a wedge about the ray in the target space on which there is a
well-defined branch of the inverse $f^{-1}_\zeta$, and if the $w_n$ are
chosen properly (in a way which depends on the angle of the wedge).  The
larger the wedge about the ray, the faster the method converges.  We now
state the exact theorem, although we shall defer the proof until
\secref{PLM-iteration}.

\begin{notation}
By a wedge $\Wedge{A,w}$, we mean the set 
   $$\{ z \st 0<|z|<2|w|,\ \arg w-A < \arg z < \arg w+A \}.$$
\end{notation}

\begin{theorem}   \Label{PLM-follows-ray}
  Suppose that the branch of the inverse $f^{-1}_{z_0}$ is analytic on 
  a wedge $\Wedge{A,w_0}$, with $0 < A \le \pi/2$, and let $h \le
  \frac{\sin A}{19}$. 
  Suppose also that $|f(z_0) - w_0| < h|w_0|/2$, and define
  	$$w_n     = (1-h)^n w_0, \qquad
          z_{n+1} = z_n - \frac{f(z_n) - w_{n+1}}{f'(z_n)}.$$
  Then $|f(z_n) - w_n| \le h|w_n|/2$ and
        $z_{n+1} \in f_{z_0}^{-1}\( \Wedge{A,w_n} \).$
\end{theorem}

\medskip
It should be noted that this theorem is a slight improvement of Smale's
Theorem~A in \Cite{Smale:bulletin85}. His proof is valid for all angles, but
is stated only for $A=\pi/12$ with $h=1/98$, and for $A=\pi/4$ with
$h=1/32$. For $A=\pi/12$ we can take take $h=1/74$, and for $A=\pi/4$,
$h=1/27$ is adequate. 

\bigskip
\subsection{Branched covers, inverse functions, and all that}
\Label{branched-covers-sec}
If $f$ is a polynomial of degree $d$, then $f:\C\maps\C$ is a branched
covering with branch points at the critical points $\theta_i$ of $f$.  If
$z$ is a regular point of $f$, that is, $f'(z) \ne 0$, then there is a
well-defined inverse function $f^{-1}_z$ so that $f^{-1}_z\(f(z)\) = z$.

In any neighborhood of a critical point of $f$, there cannot be a single
valued inverse; however, the behavior at such points is well understood.
Let $\theta$ be a critical point of multiplicity $k-1$.  Then we have
	$$f(z)-f(\theta) = (z-\theta)^k g(z),
	  {\qquad {\hbox{where}}\quad}g(\theta) \ne 0.$$
One can then define $k$ branches of the inverse which are analytic on a
small slit disk about $f(\theta)$.  We may, of course, choose any slit which
connects $f(\theta)$ to the boundary of the disk.  The reader is referred to
a complex analysis text for further details (for example, see
\Cite{Ahlfors:ComplexAnalysis}). 

\smallskip
Let $\{\zeta_j\}$ be the $d$ roots of $f$, represented with multiplicity.
If $\zeta_j$ is a simple root, denote by $f^{-1}_{\zeta_j}$ or $f^{-1}_j$
the unique branch of the inverse of $f$ which takes $0$ to $\zeta_j$.  On
the other hand, if $\zeta_j$ is a multiple root, we let $f^{-1}_{\zeta_j} =
f^{-1}_{j}$ be one of the branches of the inverse that take $0$ to
$\zeta_j$, taking care to account for all such branches exactly once.  This
can be done, since if $\zeta$ is a root of multiplicity $k \ge 2$, it is
also a critical point of multiplicity $k-1$, and so there are $k$ branches
of the inverse.

\smallskip
We now analytically continue each of the $f^{-1}_j$ to a maximal starlike
domain $\Omega_j$ in the target space; that is, we attempt to extend each
$f^{-1}_j$ along open rays from $0$.  When doing this, it is useful to think
of the target space as consisting of $d$ copies of $\C$, with a single
$f^{-1}_j$ associated to each one. 
When does the analytic continuation fail? Precisely when the inverse image
of a ray encounters a critical point of $f$. Refer to \figref{inverse-stacks}.

\begin{figure}[thb]
 \centerline{\psfig{figure=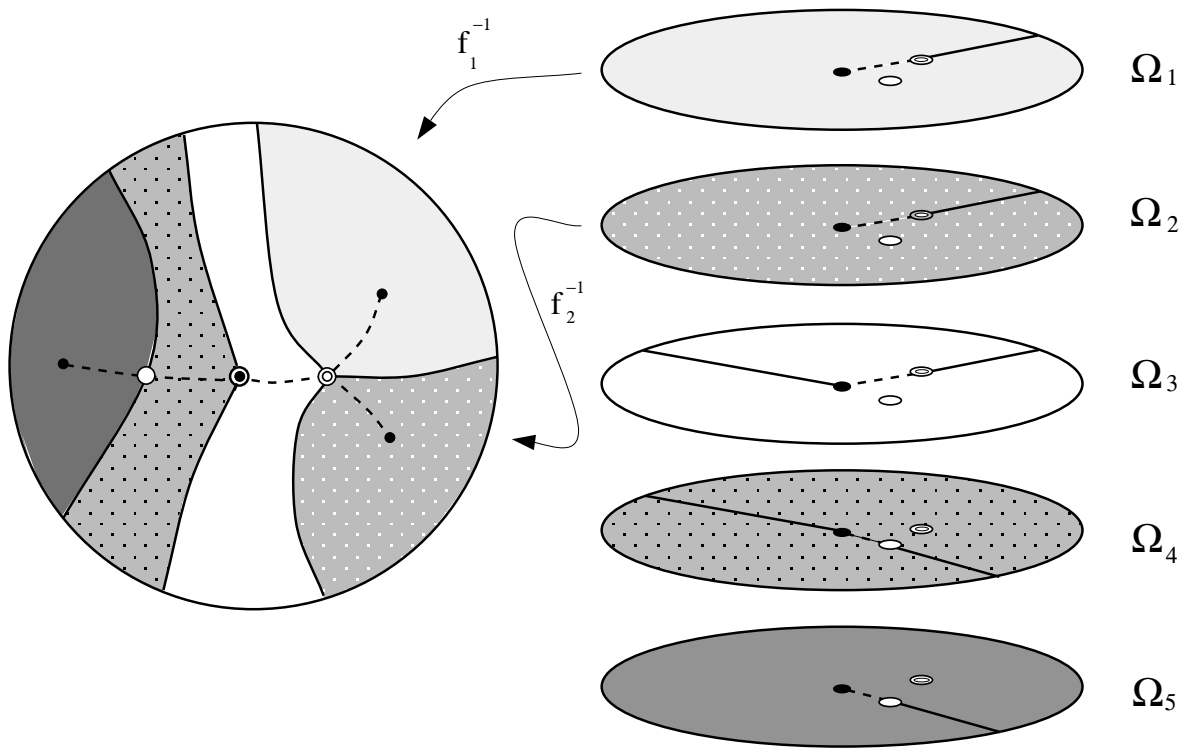,width=\textwidth}}
 \caption{\Label{inverse-stacks} The $\Omega_j$ in the target space (on the
          right), and the corresponding source space, for a degree 5
          polynomial with a double root (black dot inside a white dot)
          and a critical point of multiplicity 2 (double white dot).  The
          other critical point is marked by a single white dot, and the
          roots by black dots.  The cuts in the $\Omega_j$ are represented
          by black radial lines.}
\end{figure}

\smallskip
At this point, it may be useful to consider the Newton vector field given by
	$$\dot{z} = -{f(z)}/{f'(z)}.$$
Let $\varphi_t(z)$ be the solution curve with initial condition
$\varphi_0(z)=z$.  Then we have 
	$$f\(\varphi_t(z)\) = \e^{-t} f(z);$$
that is, $f$ maps solution curves of the Newton vector field to rays in the
target space.  Notice that the singularities of the vector field occur
precisely at the points where $f'(z)=0$.  The solution curves $\varphi_t(z)$
which have singularities play an important role here: they divide the
source space into regions on which $f$ is injective. Refer to
\Cite{Shub:PekingNotes,ShubTischlerWilliams,Smale:bulletin85} for more
details on the behavior of the solution curves.  Applying the path lifting
method can be viewed as attempting to follow the solution curves to the flow
$\varphi_t$.

\smallskip 
When constructing the $\Omega_j$, we continue $f^{-1}_j$ radially outward
until a critical point $\theta$ is encountered in the source space.  We then
exclude the ray $\{r f(\theta) \ |\ r \ge 1\}$ from $\Omega_j$, and continue
by moving along rays which avoid the cut.  Notice that when
we encounter a critical point $\theta$ of multiplicity $k-1$, we need to
slit at most $k$ of the $\Omega_j$ starting at $f(\theta).$ Also, note that
some of the $\Omega_j$ may already be slit at $f(\theta)$, since there may
be another critical point whose image lies on the same ray.

\medskip
We now count the number of such cuts:  $f$ has $d-1$ critical points (with
multiplicity), and a critical point of multiplicity $k-1$ can cause at most
$k$ cuts.  This means we have at most $2(d-1)$ cuts, distributed through the
$d$ copies of the target space.  Note that if $\Omega_j$ contains some wedge
$W$, then $f^{-1}_j$ is analytic on $W$.  The following counts the number of
$\Omega_j$  which contain wedges of a given size.

\begin{lemma}   \Label{good-sector-count}
Let $m$ be an integer, and divide $\C$ into $m$ wedges
   $$W_{n,m} = \{w \st \frac{2n\pi}{m} < \arg w < \frac{2(n+1)\pi}{m}\}
   \qquad n=0,\ldots{},m-1.$$
For each wedge $W_{n,m}$, let $N(n,m)$ be the number of $\Omega_j$ which
contain the sector $W_{n,m}$, and let 
$\displaystyle{N(m) = \max_{0\le n<m} N(n,m)}$.
Then
	$$N(m) \ge d - \left\lfloor{\frac{2(d-1)}{m}}\right\rfloor.$$
\end{lemma}

\begin{proof}
Since we have $2(d-1)$ cuts and $m$ wedges $W_n$, at least one of the wedges
has no more than $\frac{2(d-1)}{m}$ cuts.  Since there are $d$ $\Omega_j$s, 
we have the result.
\qquad\end{proof}

\begin{corollary} $N(d)=d$, $\displaystyle{N(3) \ge \frac{d}{3}}$, and 
  $\displaystyle{N(4) \ge \frac{d}{2}}$.
\end{corollary}

\begin{proof}
Application of the formula above gives the result for $N(3)$ and $N(4)$, and
yields $N(d) \ge d-1$.  However, since each critical point causes at least
two cuts, it is not possible to have a wedge cut only once.  Thus, the value
$d-1$ is not permissible for $N(d)$, giving $N(d) = d$.
\qquad\end{proof}

\medskip
\subsection{Families of root-finding methods} \Label{root-finding-families}
We now use the result of \lemref{good-sector-count} to construct families
of root-finding  algorithms.  Recall that the the modified Newton method
described in \secref{PLM-section} works when there is an $\Omega_j$
containing a wedge about our initial value $f(z_0)$; the
larger the wedge, the faster the method converges.

For each family, we start with $md$ points in the source space placed around 
a circle which contains all the roots.  We think of this as $m$ sets of $d$
initial points, and choose them so that the image of each set lies well
inside each of the $m$ sectors $W_{n,m}$ in the target space. Then by
\lemref{good-sector-count}, one of the $m$ sets of points will contain at
least $N(m)$ elements whose images are each in a  ``good wedge'', that is,
they lie in some $\Omega_j$.  As a consequence,  iterating these points
under the path lifting method will locate at least $N(m)$ roots of the
polynomial. 

Particular families of interest are $m=d$, which gives the algorithm
discussed in \Cite{Kim:TopoComplexity}, and $m=4$, on which we focus our
attention here.  The basic idea of all of the algorithms is this:  obtain
$md$ ``good'' initial points and apply the path lifting method to $d$ of
them at a time.  If, after a prescribed number of iterations, we have found
approximation to at least $N(m)$ roots (counting multiplicity), we deflate the
polynomial (that is, divide out the approximated roots) and repeat the 
process on the result.  If not, we try again with the next set of $d$ points.  
Note that we are guaranteed success by the time we try the $m^{th}$ set.

The remainder of the paper consists of a detailed description and analysis
of the algorithm for $m=4$.  Most of what follows can be readily adapted to
the other families as well.

\section{A root-finding algorithm} \label{algorithm}

\subsection{Statement of the algorithm and main theorem}
\Label{alg-statement}
Here we present our root-finding algorithm for the family $m=4$.  The
presentation is structured as a main routine and several subroutines, which 
do most of the work.

\begin{notation}
Throughout this chapter, we shall denote matrices, vectors, and sets in
uppercase calligraphic type, and their elements in subscripted lowercase
type.  For example, $x_j$ is the $j^{th}$ element of the vector ${\cal X}$. 
We shall also use the notation $\lfloor x \rfloor$ to denote the least
integer in $x$, sometimes also called {\tt floor}$(x)$.
\end{notation}

The main routine merely inputs the desired polynomial and precision,
rescales it so the roots lie in the disk of radius $1/2$, then repeatedly
calls a subroutine to halve the number of unknown roots (counted with
multiplicity) and deflate.   We do the rescaling in order to easily bound
the error introduced by the FFT deflation. The set $\Lambda$ contains all the
approximations found by the $i^{th}$ stage. 

Note that the algorithm is given for an arbitrary monic polynomial, since
only minor changes are required to normalize the input polynomial.  If it is
assumed that the input polynomial is already in $\Pd(1)$, we can take
$f_0(z) = \phi(4z)/4^d$ and  $\tau = 32\epsilon/7^{d+3}$. 

\begin{algorithm}{Main Routine}
  \stmt Input monic polynomial  ${\phi(z) = \sum_{i=0}^d a_i z^i}$
          and desired precision $\epsilon$.
  \stmt Let $f_0(z) = \phi(Kz)/K^d$, with $K=4\ \max |a_j|^{\frac{1}{d-j}}$.
  \stmt Let $\displaystyle{\tau=\frac{\epsilon}{2K^d}\(4/7\)^{d+3}}.$
  \stmt Let $i=0$ and $\Lambda=\emptyset$
  \stmt While $\#(\Lambda) < \deg(\phi)$
  \begin{itemize}
      \stmt     Let $(f_{i+1}, \Lambda_i) =$ 
            {\tt get-half-roots-and-deflate}$(f_i, \tau)$. 
      \stmt     Let $\Lambda = \Lambda \union \Lambda_i$.
      \stmt     Increment $i$.
  \end{itemize}
  \stmt End While
  \stmt Output $K \Lambda$. 
\end{algorithm}

The function {\tt get-half-roots-and-deflate} takes as input a normalized
polynomial $f$ and precision $\tau$.  It returns a set of points $y_j$ which
approximate at least half of the roots of $f$ (with multiplicity) and a new
polynomial $\tilde f$ which we obtain by deflation.  These satisfy 
$\| f(z)-\tilde{f}(z) \prod(z-y_j) \| < 2\tau$.  We should point out here
that we are actually finding approximate zeros of $f - \vec\tau$, where
$|\vec\tau| = \tau$, which depends on $\epsilon$.  When the translation is
in the proper quadrant, this will ensure that the relevant roots of
$f-\vec\tau$ are simple, so that we have approximate zeros in a
neighborhood.  This allows us to obtain the right number of approximations
to a multiple root, without worrying about winding number arguments or the
like.  We emphasize again that $\tau$ is chosen as a function of
$\epsilon$, and is small enough that the approximation polynomial has
negligible errors in the non-constant terms.  

The matrix ${\cal Z}$ consists of $4$ rows of $d$ ``good'' initial
conditions, with $|z_{j,k}|=3/2$ and $\arg f(z_{j,k}) \approx 2\pi\i j / 4$.
We use  ${\cal Z}_j$ to represent the $j^{\rm th}$ row. 

\begin{algorithm}{function get-half-roots-and-deflate$(f, \tau)$}
  \stmt Let $d = \deg f$.
  \stmt Let ${\cal Z} =${\tt choose-4d-good-initial-points}$(f)$. 
  \stmt For $j=1$ to $4$ do
      \begin{itemize}
      \stmt Let ${\cal Y} =$ {\tt iterate-PLM}$(f,{\cal Z}_j,j,\tau)$.
      \stmt Let $\psi = f-\tau \e^{2\pi\i j /4}$.
      \stmt Let ${\cal X} =$ {\tt select-approx-zeros}$(\psi,{\cal Y})$.
      \stmt Let ${\cal W} =$ {\tt polish-roots}$(\psi,{\cal X},\tau)$.
      \stmt Let ${\cal V} =$ {\tt weed-out-duplicates}$(\psi,{\cal W})$.
      \stmt If $\#({\cal V}) \ge d/2$ then
      \begin{itemize}
        \stmt Let $\tilde f=$ {\tt deflate}$(\psi,{\cal V})$.
        \stmt Return ($\tilde{f}, {\cal V}$)
      \end{itemize}
      \stmt End if
    \end{itemize}
  \stmt End for
  \stmt Print ``We have proven this statement will not be reached.  Check
         your code.''
  \stmt Abort.
\end{algorithm}

The function {\tt choose-4d-good-initial-points} gives us our 4 sets of $d$
initial values on the circle of radius $3/2$.  The sets have the property
that elements of the same set are mapped to points with approximately the
same argument, and elements of distinct sets are mapped to points in
different quadrants of the target space.  There are several different ways
to accomplish this, and other methods can slightly decrease the number of
operations required in {\tt iterate-PLM}.  Please refer to the remarks in 
\secref{initial-points} for more details.

\begin{algorithm}{{\tt function choose-4d-good-initial-points$(f)$}}
  \stmt Let $N = 676 \deg f$.
  \stmt For $k=1$ to $N$
  \begin{itemize}
    \stmt Let $\omega_k = \frac{3}{2} \e^{2 \pi\i k/N}$.
  \end{itemize}
  \stmt End For.
  \stmt For $j=1$ to $4$
  \begin{itemize}
    \stmt Let ${\cal Z}_j$ be the union of the $\omega_k$ for
    which $\arg f(\omega_k) \le 2\pi/j$ and $\arg f(\omega_{k+1}) > 2\pi/j$.
  \end{itemize}
  \stmt End For.
  \stmt Return(${\cal Z}$).
\end{algorithm}

From our sets of initial points ${\cal Z}_j$, we obtain our approximate
zeros via the routine {\tt iterate-PLM}.  This applies the path lifting
method to $f$ for an appropriate number of steps.  For simplicity, we present
a scalar version here.  However, there are well-known methods for
evaluating the same polynomial $f$ at $m$ different points ($m \ge d$) in
$\Order{m(\log d)^2}$ arithmetic operations; refer to section~4.5 of
\Cite{BorodinMunro}. 
When computing the complexity in \secref{complexity}, we assume such methods
are used.  Note that this algorithm can be easily implemented on either a
vector or parallel computer, so that one can iterate the $d$ elements of
${\cal Z}_j$ simultaneously.  Also note that if Smale's mean value
conjecture holds, the extra iteration to obtain $\hat{z}$ is not needed;
$z_N$ is good enough.  See the remarks following
\lemref{approx-zeros-if-radius-small}. 

\begin{algorithm}{{\tt function scalar-iterate-PLM$(f,z_0,j,\tau)$}}
  \stmt Let $w_0 = |f(z_0)| \e^{j \pi \i/2}$.
  \stmt Let $h = 1/27$.
  \stmt Let $\displaystyle{N =\left\lfloor{
           \frac{\log_2(\tau/|w_0|)}{\log_2(26/27)} 
         }\right\rfloor}$
  \stmt For $n = 1$ to $N$
  \begin{itemize}
        \stmt Let $w_n = (1-h)w_{n-1}$.
        \stmt Let $\displaystyle{z_n = z_{n-1} - 
                                \frac{f(z_{n-1})-w_i}{f'(z_{n-1})}}$.
  \end{itemize}
  \stmt End for.
  \stmt Let $\displaystyle{\hat{z} = z_N - 
                           \frac{f(z_N)-\tau \e^{j \pi \i/2}}{f'(z_N)}}$.
  \stmt Return$(\hat{z})$.
\end{algorithm}

The next routine takes the output of {\tt iterate-PLM} and uses the
$\alpha$ function (defined in \lemref{alpha}) to remove those
elements which are not approximate zeros.  Although the test $\alpha < 1/8$
is sufficient, it is not a necessary condition.  However, if the image of
the initial points ${\cal Z}_j$ lie in a ``good quadrant'', we know by 
\lemref{PLM-its} below that we will have $\alpha < 1/8$ for at least half of
them, and so we are not in danger of discarding too many points.

\begin{algorithm}{{\tt function select-approx-zeros$(\psi,{\cal Y})$}}
  \stmt Let ${\cal X}=\emptyset$.
  \stmt For $j=1$ to $\#({\cal Y})$
  \begin{itemize}
    \stmt If $\alpha_\psi(y_j) < 1/8$, then let ${\cal X} = {\cal X} \union
          \{y_j\}$. 
  \end{itemize}
  \stmt End For.
  \stmt Return(${\cal X}$).
\end{algorithm}

Once we have found approximate zeros for at least half of the roots of
$f_i$ (by applying {\tt iterate-PLM} to at most 4 sets ${\cal Z}_j$), we
further refine them by ``polishing'' with regular Newton's method.  As in {\tt
iterate-PLM}, fast polynomial evaluation techniques can be used, but we
present a scalar version here for simplicity.

\begin{algorithm}{{\tt function scalar-polish-roots$(\psi,x_0,\tau)$}}
  \stmt Let $\displaystyle{M = 1 + \left\lfloor{\log_2 \log_2 
        \(64d(7/4)^d\)/{\tau} }\right\rfloor}.$       
  \stmt For $n=1$ to $M$
  \begin{itemize}
    \stmt Let $\displaystyle{x_n = x_{n-1} -
          \frac{\psi(x_{n-1})}{\psi'(x_{n-1})}}$.
  \end{itemize}
  \stmt End For.
  \stmt Return$(x_M)$.
\end{algorithm}

Finally, we remove from the set of approximations that we have found any
points which approximate the same root of $\psi$.  We do this by taking
the approximation which minimizes $|\psi|$, and then making a pass through
the rest of them and accepting only those which approximate different roots
from the ones previously accepted. 

\begin{algorithm}{{\tt function weed-out-duplicates$(\psi,{\cal W})$}}
  \stmt Sort ${\cal W}$ so that 
        $|\psi(w_1)| \le |\psi(w_2)| \le \ldots{} \le |\psi(w_n)|.$
  \stmt Let ${\cal V}=\{w_1\}$.
  \stmt For $j=2$ to $\#({\cal W})$
  \begin{itemize}
      \stmt If $|w_j - v| > 3|\psi(w_j)|/|\psi'(w_j)|$ for all 
            $v \in {\cal V}$, then let ${\cal V} = {\cal V} \union \{w_j\}$. 
  \end{itemize}
  \stmt End For.
  \stmt Return(${\cal V}$).
\end{algorithm}

\medskip
At this point we have found approximations to at least half of the roots of
$f_i$. We divide them out to obtain a new polynomial $f_{i+1}$ of smaller
degree, using a standard technique involving the finite Fourier matrix.

\begin{algorithm}{{\tt function deflate($\psi,{\cal V}$)}}
  \stmt Let $n = \deg \psi - \#{\cal V}$.
  \stmt Let $\omega = \e^{2\pi\i /(n+1)}$.
  \stmt Let $\displaystyle{p(z) = \prod_{k=1}^{\#{\cal V}}\(z - v_k\)}$.
  \stmt Let ${\cal M}^{-1}$ be the inverse of the $n\times n$ Fourier
        matrix.  That is, $m_{j,k} = \omega^{-j k}/n$.
  \stmt Let ${\cal P}$ be the column vector with entries
        $\displaystyle{p_j=\frac{\psi(\omega^j)}{p(\omega^j)}},
        \ \ j=0,\ldots{},n$.
  \stmt Let ${\cal Q} = {\cal M}^{-1}\ {\cal P}$.
  \stmt Let $q(z) = \sum_{j=0}^{n} q_j \, z^j$.
  \stmt Return($q$).
\end{algorithm}

\medskip
We now state our main theorem, which essentially says that the algorithm
just presented works:

\begin{theorem} \Label{main-theorem}
  Given a monic polynomial $\phi(z)$ of degree $d$ and $\epsilon > 0$, the
  algorithm presented in this section will always terminate 
  with $d$ points $\lambda_1,\ldots{},\lambda_d$ which satisfy
  $$\| \phi(z) - \prod_{j=1}^{d}(z - \lambda_j) \| < \epsilon.$$
  For $\phi \in \Pd(1)$, the worst case arithmetic complexity of this
  algorithm is 
  	$$\Order{d(\log d)^2|\log\epsilon| + d^2(\log d)^2}.$$ 
\end{theorem}

\begin{corollary}
  Let $\phi$ and $\epsilon$ be as in the theorem, and denote the roots of
  $\phi$ by $\xi_j$, represented with multiplicity.  Then the algorithm can
  be used   to produce the approximations $\lambda_j$ satisfying
  	$$| \xi_j - \lambda_j | < \epsilon, $$
  with a worst-case arithmetic complexity of
  	$$\Order{d^2(\log d)^2|\log\epsilon| + d^2(\log d)^2}.$$ 
\end{corollary}

\begin{proof}
This is an immediate consequence of the fact that if $f$ and $g$ are in
$\Pd(1)$ with  $\| f - g \| < \(\epsilon/8d\)^d$, then the roots of $f$ and
$g$ are at most $\epsilon$ apart (See \cite{Kim:TopoComplexity}).
\qquad\end{proof}

We prove the theorem as a series of lemmas in the following sections.  There
is a rough correspondence between the sections and the routines in the
algorithm.  Finally, we summarize all of these lemmas and give the proof in
\secref{proof-of-main-thm}. 

\medskip
\subsection{Selection of initial points}        \Label{initial-points}
For each intermediate polynomial $f_i$, we need to select four sets of
$(\deg~f_i)$ points at which to begin our iteration.  These must be chosen
so that the elements of each set map very near the same point in the target
space, and that the images of elements in successive groups are approximately
$\frac{1}{4}$-turn apart.  This can be accomplished either by evaluating
$f_i$ at a large number of points spaced around a circle in the source space,
and then selecting from those, or by taking a much smaller number (perhaps
as few as $4d$) of points and adjusting them with either a standard or
modified  Newton's method.   
The arithmetic complexity of either comes out much the same; we
opt for the former method because of its conceptual simplicity.

The following lemma gives us bounds on how much the argument in the target
space can vary between points around a circle containing all the roots in
the source space.  Using this, we see how many points are required to
obtain our ``good'' points. 

\begin{lemma}   \Label{init-arg}
  Suppose all of the roots $\zeta$ of a polynomial $f$ of degree $d$ lie in
  $\D_R(0)$, and let 
	$$\omega_m = 2 R\, \e^{\frac{2 \pi \i m}{n d}}.$$
  Then
  $$ \left| \arg \frac{f(\omega_{m+1})}{f(\omega_{m})} \right | 
     \le \frac{4 \pi}{n}.  $$
\end{lemma}

\begin{proof}
  This argument appears in \Cite{Renegar}(Lemma 7.1), although in a somewhat
  different form.  We present an adapted version here.  The idea is quite
  simple: since 
$$      \arg \frac{f(\omega_{m+1})}{f(\omega_{m})} =
	\arg \frac{\prod_{i=1}^d (\omega_{m+1}-\zeta_i )}
        	  {\prod_{i=1}^d (\omega_{m}-\zeta_i)} =
        \sum_{i=1}^d \arg  \frac{\omega_{m+1}-\zeta_i}{\omega_{m}-\zeta_i}, $$
  we merely bound the angles in the source space and add them up.

\begin{figure}[thb]
 \centerline{\psfig{figure=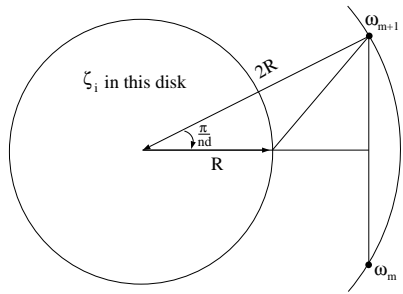,height=1.5in}}
 \caption{\Label{init-arg-fig} Calculation of the upper bound on 
 the variation of argument in the target space.  This picture is in the
 source space.}
\end{figure}

  Note that each root $|\zeta_i|<R$ and $|\omega_m|=2R$, so
   \begin{eqnarray*}
          \left| \arg  \frac{\omega_{m+1}-\zeta_i}{\omega_{m}-\zeta_i} \right |
    &\le& 2 \arctan \frac { |\omega_{m+1}-\omega_m| / 2 } 
    			  {2R \cos\frac{\pi}{nd} - R } \\
    &\le& 2 \arctan \frac { 2 R \sin \frac {\pi}{nd} }
    			  {2 R ( \cos\frac{\pi}{nd} - \frac{1}{2}) } \\
    &\le& \frac {4\pi}{nd}.
   \end{eqnarray*}
  Refer to  \figref{init-arg-fig}. Summing the $d$ terms gives the result.
\qquad\end{proof}

From \thmref{PLM-follows-ray}, for $h=1/27$ we require that our initial
points be within $\sin^{-1}(h/2) < \pi/169$ of the central ray, so we start
with $676d$ points equally spaced around the circle of radius $3/2$.  We then
evaluate the polynomial at each of them, and select four sets of $d$ points
whose arguments are closest to $0$, $\frac{\pi}{2}$, $\pi$, and
$\frac{3\pi}{2}$, respectively.  For each of these, we take the initial
target point $w_0$ to be the projection of its image onto the real or
imaginary axis. 

\begin{remark}
Some amount of computation can be saved if we use the same $w_0$ for all $d$
elements of a given set of initial points, rather than just points with the
same argument.  This would make the computation of the $w_n$ a scalar rather
than a vector operation, that is, $w_n$ would only need to be computed once
for each group of $d$ points.  However, in order to do this, we must ensure
that the images of each $z_j$ in the same set have approximately the same
modulus as well as argument.  This is best accomplished using some sort of
Newton's method. 

If an initial Newton's method is used, one can also choose a much smaller
number of trial points $\omega_i$.  Such a method should converge well, since
all the critical points of $f(z) - w_0$ are inside the $\D_R$, while the
roots of $f(z)-w_0$ and the $\omega_i$ are well outside.
\end{remark}

\medskip
\subsection{Iteration of the path lifting method} \Label{PLM-iteration}
In this section, we analyze the behavior of applying the path lifting
method to a single well-chosen initial value $z_0$. 

We first prove the theorem as promised in \secref{PLM-section}, and then we
show that after a specified number of iterations, the result will be an
approximate zero.  Before doing this, we shall state the relevant special
case of Theorem~3.2 of \Cite{Kim:ApproxZeros}, which gives a lower bound on
how far a point moves under Newton's method.

\begin{lemma} \Label{B-lemma}
Let $\hat{z} = z - g(z)/g'(z)$, with $r = |g(z)|/R_g(z) < 0.148$, where
$R_g(z)$ is the radius of convergence of $g^{-1}_z$ as a power series based
at $z$.  Then 
$$ \left| \frac{g(\hat{z})}{g(z)} \right| \le B(r) 
   \qquad {\hbox{where}} \qquad
   B(r) = 2r\frac{(1+r)^3}{(1-r)^5}.$$
\end{lemma}

Given this lemma, the proof of \thmref{PLM-follows-ray} is quite simple.
Note that each step of the iteration in  \thmref{PLM-follows-ray}
corresponds to Newton's method applied to $f(z) - w_{n+1}$ at $z_{n-1}$.

\begin{plmthm}
  Suppose that the branch of the inverse $f^{-1}_{z_0}$ is analytic on 
  a wedge $\Wedge{A,w_0}$, with $0 \le A \le \pi/2$, and let $h \le
  \frac{\sin A}{19}$. 
  Suppose also that $|f(z_0) - w_0| < h|w_0|/2$, and define
  	$$w_n     = (1-h)^n w_0, \qquad
          z_{n+1} = z_n - \frac{f(z_n) - w_{n+1}}{f'(z_n)}.$$
  Then $|f(z_n) - w_n| \le h|w_n|/2$ and
        $z_{n+1} \in f_{z_0}^{-1}\( \Wedge{A,w_n} \).$
\end{plmthm}

\begin{proof}
  We shall prove this by induction.  All that is required is to establish
  the conclusion, given that $|f(z_{n-1}) - w_{n-1}| < h|w_{n-1}|/2$.  

  Note that
    $$    |f(z_{n-1}) - w_n| \le |f(z_{n-1}) - w_{n-1}| + |w_{n-1} - w_n|
    \le   3 h |w_{n-1}| / 2, $$
  and that
    $$    R_{f-w_n}(z_{n-1}) = R_f(z_{n-1}) \ge |w_n|\sin A - h|w_n|/2. $$
  Since $h \le \frac{\sin A}{19}$, we can apply \lemref{B-lemma} with 
  $g(z)=f(z)-w_n$ and $r=3/37$ to obtain
  $$
    \left| \frac{f(z_n) - w_n}{w_n} \right| =
    \left| \frac{f(z_n) - w_n}{f(z_{n-1}) - w_{n}} \right| \cdot
    \left| \frac{f(z_{n-1}) - w_{n}} {w_{n-1} (1-h)} \right|  
    \le    B(r) \frac{3h}{2(1-h)}.
  $$
  Because $B(3/37) < (1-h)/3$, we have our  conclusion.
\qquad\end{proof}

\begin{remark}
  The denominator of 19 in the upper bound on $h$ can be relaxed only
  slightly in this proof. If we take $h \le \frac{\sin A}{k(A)}$,
  then $k(A)$ is a monotonically increasing function, with $18.3096 \le k(A)
  \le 18.895$.
  Also, if $A > \pi/2$, we can take $h = 1/19$.
\end{remark}

We apply the path lifting method until we have an approximate zero, at
which time we can use the standard Newton iteration.  Since we are
interested in finding $\tau$-roots of $f$,
we stop the iteration once
we have an approximate zero for $f-\tau\frac{w_0}{|w_0|}$.  Note that, as
discussed in \secref{alg-statement}, this translation is necessary, since
there are no approximate zeros in a neighborhood of a multiple root.  Having
approximate zeros for the perturbed polynomial gives us the following bound
on the number of iterations required, and enables the weeding out of
duplicates. 

\begin{lemma}      \Label{PLM-its}
  Let $h\le \frac{\sin A}{19}$,  with $w_i$ and $z_i$ as in
  \thmref{PLM-follows-ray}.  If  
    $$N =\left\lfloor{
           \frac{\log_2(\tau/|w_0|)}{\log_2(1-h)} 
         }\right\rfloor, $$
  then $z_N$ is an approximate zero for $f - \tau\frac{w_0}{|w_0|}$.
  Furthermore, if
  $\hat{z} = z_N - \frac{f(z_N)-\tau\frac{w_0}{|w_0|}}{f'(z_N)}$
  then we have 
  $\alpha(\hat{z}) \le 1/8$.
\end{lemma}

\begin{proof}
  For notational convenience, set $\vec\tau = \tau\frac{w_0}{|w_0|}$.

  By assumption,  $f^{-1}_{z_0}$ is   analytic on the wedge
  $\Wedge{A,w_0}$, and the initial point 
  $z_0$ is close enough to $w_0$ that we can apply \thmref{PLM-follows-ray}.
  Note that if $N$ is as specified, we have 
  	$$ \tau \le |w_N| \le \frac{\tau}{1-h}, $$
  since $|w_N| = (1-h)^N|w_0|$ by definition.  Also, $|f(z_N) - w_N| \le
  h|w_N|/2$.  Refer to \figref{PLM-its-fig}.

\begin{figure}[thb]
 \centerline{\psfig{figure=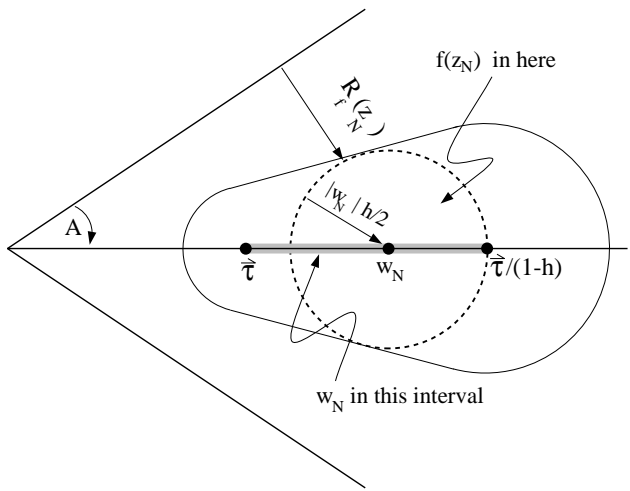,height=2.5in}}
 \caption{\Label{PLM-its-fig}The locations of $w_N$ and $f(z_N)$ in the
 target space when $z_0$ is in a ``good quadrant''.}
\end{figure}

In order to conclude that $z_N$ is an approximate zero, we need to show that
we have $\left|f(z_N) - \vec\tau\right|/R_f({z_N}) < 1/10$, where
$R_f({z_N})$ is the radius of convergence of $f^{-1}_{z_N}$. 
Consider $f(z_N)$ in some circle of radius $h|w|/2$, where $\arg w = \arg
w_0$,  and $|w| > \tau$.  Then we have  
$$
      \frac{\left|f(z_N) - \vec\tau\right|}{R_f({z_N})} \le
      \frac{|w| + h|w|/2 - \tau }{|w|\sin A - h|w|/2}             =
      \frac{2 + h}{2\sin A - h} - \frac{\tau}{|w|(\sin A - h/2)}. 
$$
The maximum of this quantity occurs when $|w|$ is as large as possible,
which in our case is $|w| = \tau/(1-h)$.  Hence
$$
        \frac{\left|f(z_N) - \vec\tau\right|}{R_f({z_N})} \le
        \frac{3h}{2\sin A -h} \le 3/37 < 1/10.$$
By \lemref{approx-zeros-if-radius-small}, $z_N$ is an approximate zero of 
$f(z) - \vec\tau$.

In order to ensure that $\alpha < 1/8$, we need to apply Newton's method
once. Set
$$      \hat{z} = z_N - \frac{f(z_N)-\vec\tau}{f'(z_N)}.$$
Then this point satisfies
$$ 
      \frac{\left|f(\hat{z}) - \vec\tau\right|}
           {\left|f(z_N) - \vec\tau\right|}        \le
      B(3/37) \le 0.31271, $$
and so 
$$      \left|f(\hat{z}) - \vec\tau \right|    \le
	0.31271 \( |f(z_N) - w_N| + |w_N| - \tau \)           \le
        0.31271 \frac{3h\tau}{2(1-h)} \le \frac{\tau\sin A}{39}.      $$
This gives us
$$ \frac{\left|f(\hat{z}) - \vec\tau\right|}{R_f({\hat{z}})} \le
	\frac{\tau\sin A / 39} {\tau\sin A - \tau\sin A /39}        = 1/38. $$
As a consequence of \lemref{approx-zeros-if-radius-small}, we have
$\alpha(\hat z) < 1/8$. 
\qquad\end{proof}

\medskip
\subsection{Refinement of the root approximations}
\Label{weeding}
The routine {\tt iterate-PLM} outputs a set of points ${\cal Y} = 
\{y_1,y_2, \ldots{}, y_d\}$ which may be approximate zeros for
$\psi(z) = f(z)-\vec\tau$.  We can use $\alpha_\psi$ (see \lemref{alpha})
to choose those $y_i$ that are indeed approximate zeros;  we discard those
$y_i$ for which $\alpha_\psi(z) > 1/8$.  As a consequence of
\lemref{approx-zeros-if-radius-small} and \lemref{PLM-its}, those $y_i$
which started  in the ``good sector'' will not be discarded. 

In addition, we also want to ensure that we approximate each root only once
(counted with multiplicity). \lemref{distinct-zeros} 
gives us conditions which allow us to weed out any duplicates; the proof
relies on a variant of the Koebe Distortion Theorem which we quote here from
\Cite{Kim:ApproxZeros}, Lemma~3.3. 

\begin{lemma}  \Label{Koebe-distortion}
  Let $g$ be univalent on $\D_R(z)$.  Then for any $s<1$, we have
     $$\D_{tR}\(g(z)\) \subset g\(\D_{sR}(z)\) \subset \D_{uR}\(g(z)\),$$
  where $t=s|g'(z)|/(1+s)^2$ and $u=s|g'(z)|/(1-s)^2.$
\end{lemma}

\medskip
\begin{lemma}   \Label{distinct-zeros}
  Suppose that $y_1$ and $y_2$ are approximate zeros of $\psi(z)$, with
  $|\psi(y_1)| \ge |\psi(y_2)|$, and $R_\psi(y_j) \ge 12\psi(y_j)$, where
  $R_\psi(y_j)$ is the radius of convergence of $\psi^{-1}_{y_j}$.  Then
  $y_1$ and $y_2$ approximate the same simple root $\xi$ of $\psi$ if and
  only if
$$
	|y_1 - y_2| < \left| \frac{3\psi(y_1)}{\psi'(y_1)} \right|.
$$
\end{lemma}

\begin{proof}
  If $y_1$ and $y_2$ approximate the same simple root, then we have
  $\psi^{-1}_{y_1} = \psi^{-1}_{y_2}$. Since $|\psi(y_1) - \psi(y_2)| \le
  2|\psi(y_1)|$,  we have $\psi(y_2) \in \D_{2|\psi(y_1)|}(\psi(y_1))$.  Thus
  we can apply \lemref{Koebe-distortion} with $g=\psi^{-1}$, which is
  univalent on the disk of radius $R=12|\psi(y_1)|$.  Taking $s=1/6$, we
  have   $\psi^{-1}\(\D_{2|\psi(y_1)|}(\psi(y_1))\)$ contained in the disk of 
  radius  $\frac{72|\psi(y_1)|}{25|\psi'(y_1)|}$ about $y_1$, so 
  $|y_1 - y_2| < \frac{3|\psi(y_1)|}{|\psi'(y_1)|}$. 

  For the other direction, we apply the Koebe $\frac{1}{4}$-Theorem (or
  \lemref{Koebe-distortion} with $s=1$; it is the same thing) to see that
  $\psi^{-1}_{y_1}\(\D_{12|\psi(y_1)|}(\psi(y_1))\)$ contains the disk of
  radius ${3|\psi(y_1)|}/{|\psi'(y_1)|}$ about $y_1$.  Thus, if the
  distance between  $y_1$ and $y_2$ is less than this amount, we must have
  $\psi^{-1}_{y_1}(\psi(y_2)) = y_2$, that is, $y_1$ and $y_2$ approximate 
  the same root of $\psi$. 
\qquad\end{proof}

\begin{remark}
Note that $\alpha_{\psi}(z) < 1/8$ is not sufficient to imply that the
hypotheses of \lemref{distinct-zeros} are satisfied, since this only gives
$R_\psi(z) \ge 4|\psi(z)|/3$.  However, if Newton's method is applied to such
a point at least 3 times, the value of $|\psi|$ will decrease by at least
$1/128$, and so the lemma can be applied.   Since we need to ``polish'' the
approximations with Newton's method in order to control the error in the
deflation, we do that before weeding out the duplicates.  The total number of
iterations of Newton's method required is calculated in \secref{err-bound},
but it is greater than 3 in all cases.  In practice, one should 
probably perform the weeding in the routine {\tt select-approx-zeros}, in
order to avoid polishing points which will be discarded later.
\end{remark}

\medskip
\subsection{Deflation of intermediate polynomials}
Here we compute an explicit bound on the error introduced by the deflation
step.  We start with a polynomial $\psi$ which has roots
$\{\xi_j\}_{j=1,\ldots{},d}$, and a set of approximations to these roots which
we denote by  $\{v_k\}_{k=1,\ldots{},n}$ with $n\le d$.  We then use polynomial
interpolation via the discrete Fourier matrix to obtain a polynomial $q$ of
degree $d-n$ so that 
	$$\psi(z) \approx p(z) q(z),\qquad 
	\hbox{where}\qquad  p(z) = \prod_{k=1}^n (z - v_k).$$
In this section, we estimate $\| \psi(z) - p(z) q(z) \|$, as well as
the accumulated error in repeating this process until $q(z)$ is a constant.
\smallskip
Recall that the  finite Fourier matrix ${\cal M}$ has as its
$j,k^{th}$ entry the $jk^{th}$ power of a primitive $d+1^{st}$-root of unity.
That is,
	$$m_{j,k} = \omega^{jk} = \e^{2\pi\i jk/(d+1)}
        \qquad j,k=0,\ldots{},d$$
Then one can readily see that if $f(z) = \sum_{j=0}^d a_j z^j$ and
${\cal A}$ is the column vector of the coefficients of $f$, then the product 
${\cal M \cal A}$ will be the vector ${\cal B}$ whose $j^{th}$ entry is the
value of $f(\omega^j)$.  Also, given the values of $f$ evaluated at the
powers of $\omega$, we can easily compute the coefficients of $f$ as the
product  ${\cal M}^{-1}{\cal B}$.

\smallskip
Our first lemma gives us estimates on the size of the error caused by a single
deflation.
As is common with this sort of thing, the proof is neither terribly
entertaining or enlightening. 

\begin{lemma} \label{single-deflation}
Suppose $\psi(z) = \prod_{j=1}^d (z - \xi_j)$, with $|\xi_j| < 3/4$, factors
as $\psi(z) = P(z)Q(z)$, where $\deg(P) = n < d$ and $\deg(Q)=m=d-n$.
Let $p(z)= \prod_{j=1}^n ( z - v_j)$, with 
$| \xi_j - v_j| < \delta < \frac{1}{8d^2}$, and define $q(z)$ and $r(z)$ by
$\psi(z) = p(z)q(z) + r(z)$, where $q$ is found by polynomial interpolation
as described above.  Then
   \begin{eqnarray*}
    \|P(z) - p(z)\|     <\, & 8n\delta\(7/4\)^n    \\
    \|Q(z) - q(z)\|     <\, & 8m\delta\(7/4\)^m    \\
    \|r(z)\|            <\, & 8d\delta\(7/4\)^d    .
   \end{eqnarray*}
\end{lemma}

\begin{proof}
  Let $\omega = \e^{2\pi\i/(n+1)}$,  and let ${\cal B}$ be the vector with
  entries $b_j = P(\omega^j) - p(\omega^j)$, for $0 \le j \le n$.  Note that
  \begin{eqnarray*}
    b_j &=& \prod_{k=1}^n (\omega^j - \xi_k) - \prod_{k=1}^n (\omega^j - v_k) \\
        &=& \( \prod_{k=1}^n (\omega^j - \xi_k) \)                             
           \( 1 - \prod_{k=1}^n \frac{\omega^j - v_k}{\omega^j - \xi_k}  \)   \\
        &=& \( \prod_{k=1}^n(\omega^j - \xi_k) \)
           \( 1 - \prod_{k=1}^n\(1+ \frac{\xi_k - v_k}{\omega^j - \xi_k} \)\).
  \end{eqnarray*}
  Also note that 
  $$\frac{1}{4} < |\omega^j - \xi_k| < \frac{7}{4}
    \qquad {\hbox{and}} \qquad
    \left| \frac{\xi_k - v_k}{\omega^j - \xi_k} \right| \le 4\delta .$$
  Since $4\delta \le \frac{1}{2d^2} < \frac{1}{2n^2}$, we have
  $(1+4\delta)^n \le 8n\delta +1$,  and so
  $$|b_j| < \(7/4\)^n \bigr( (1+4\delta)^n - 1\bigl) 
  	  < 8n\delta\(7/4\)^n.$$
  Thus
  $$ \|P(z)-p(z)\| = \|{\cal M}^{-1}{\cal B}\| \le \|{\cal B}\|
  	           < 8n\delta\(7/4\)^n.$$

\bigskip
  We have a similar computation for the bound on $\|Q-q\|$:
  Let $\eta = \e^{2\pi\i/(m+1)}$ and let ${\cal C}$ be the vector with
  entries  
  $c_j = Q(\eta^j) - q(\eta^j),\ j=0,\ldots{},m$.  
  Then
  $$ c_j = \frac{\psi(\eta^j)}{P(\eta^j)} - \frac{\psi(\eta^j)}{p(\eta^j)} 
  	 = \frac{\psi(\eta^j)}{P(\eta^j)}  \(1 - \frac{P(\eta^j)}{p(\eta^j)} \)
         = Q(\eta^j) \(1 - \prod_{k=1}^m\frac{\eta^j - v_k}{\eta^j - \xi_k}\).$$
  As in the case of $P(z) - p(z)$, we obtain 
  $$
	\|Q(z)-q(z)\| = \|{\cal M}^{-1}{\cal C}\| \le \|{\cal C}\|
  	           < 8m\delta\(7/4\)^m.
  $$

  \bigskip
  Finally, for the bound on $\|r(z)\|$, note that $r(z) = P(z)Q(z) -
  p(z)q(z)$.  If we write $r(z) = \sum_{j=0}^d r_j z^j$, $P(z) =
  \sum_{j=0}^d P_j z^j$, 
  and so on, then we have
  $$
   r_j = \sum_{k=0}^j P_{j-k}Q_k - \sum_{k=0}^{j} p_{j-k}q_k
       = \sum_{k=0}^j \bigr( q_k (P_{j-k} - p_{j-k}) - P_{j-k} (Q_k - q_k)\bigl).
  $$
  Since $|\xi_j| < 3/4$, we have the following crude bounds on the
  coefficients of $P$ and $q$: 
  $$ |P_{j-k}| \le \(7/4\)^n
     \qquad {\hbox{and}}\qquad
     |q_j| \le  \(7/4\)^m.
  $$
  Combining this with the bounds on $\|P-p\|$ and $\|Q-q\|$, we get
    \begin{eqnarray*}
     |r_j| &\le&\ \sum_{k=0}^j \( 8m\delta\(7/4\)^m \(7/4\)^n
     			    +  8n\delta\(7/4\)^n \(7/4\)^m\)    \\
           &\le&\ 8j d \delta \(7/4\)^d                          \\
           &\le&\ 8d^2  \delta \(7/4\)^d.
    \end{eqnarray*}
\end{proof}

Now that we have bounds on the error in one step of deflation, we can
bound the error introduced by repeated deflation.  We assume that our
initial polynomial $f$ has roots in $\D_{1/2}$ so that we can ensure
that the roots of the subsequent polynomials $f_k$ remain in
$\D_{3/4}$ as required by  \lemref{single-deflation}. 

\begin{lemma} \Label{repeated-deflation}
Suppose $f(z) = \prod_{j=1}^d (z - \zeta_j)$, with $|\zeta_j| < 1/2$.
Let $f_0 = f$, and define 
$$      f_k(z) = p_{k+1}(z) f_{k+1}(z) + r_{k+1}(z) $$
for $0<k<m-1$, where $p_j$, $f_j$, and $r_j$ are determined by polynomial
interpolation as in \lemref{single-deflation}, with $f_m(z) = 1$.  Suppose
also that $\deg r_{k+1} < \deg f_k=d_k$, and that 
$\|r_k\| < \min\(\frac{1}{4(4d)^d},\mu({4/7})^{d+3}\).$
Then
$$      \| f - p_1 p_2 \cdots{} p_m \| < \mu.$$
Furthermore, if we have  $d_{k} \le d_{k-1}/2$ and $\mu \le (7/16)^d$, we
need only require that $\|r_k\| < \mu(4/7)^{d+3}$.
\end{lemma}

\begin{proof}
Note that $f-p_1 p_2 \cdots p_m = p_1\cdots p_m r_{m-1} + \ldots{} + p_1 r_2
+ r_1 $,  and so
$$\|f-p_1 p_2 \cdots p_m\| \le \| p_1\cdots p_{m-1}\|\|r_{m-1}\| +
  \ldots{} + \|p_1\|\|r_2\| + \|r_1\|.$$
First, we show that $\|p_1\cdots{}p_k\| \le (7/4)^{\deg p_1\cdots{}p_k}$.
We shall use induction to show that the roots $p_1\cdots{}p_k$ are
always in the circle of radius $\frac{1}{2} + \frac{k}{4d} < 3/4$.
On the circle $|z| = \frac{1}{2} + \frac{k}{4d}$,  we have
$$ |f_{k-1}(z) - p_k(z) f_k(z)| 
   \le \|r_k\| \sum_{j=0}^{\deg r_k}|z|^j
   \le \frac{d}{(4d)^d (2d-k)} 
   \le \frac{1}{(4d)^d}  \le |f_{k-1}(z)|, $$
and so by Rouch\'e's Theorem
the roots of $p_k f_k$ lie inside $\D_{\frac{1}{2} + \frac{k}{4d}}$.  
Thus the roots of $p_1 p_2 \cdots{} p_k$ also lie in this disk, and so the
coefficients are less than $(7/4)^{\deg p_1 \cdots{} p_k}$.

\smallskip
Since $\deg p_1 p_2 \cdots{} p_k < d-(m-k)$, we can conclude that
\begin{eqnarray*}
  \sum_{k=1}^{m-1} \|p_1 p_2 \cdots{} p_k\| \|r_{k+1}\| 
  &\le&\ \sum_{k=1}^{m-1} (7/4)^{d-(m-k)} \(4/7\)^{d+3}\mu  \\
  &\le&\ \frac{4}{3}  \(7/4\)^{d+2} \(4/7\)^{d+3}\mu        \\
  &<&\ \mu.
\end{eqnarray*}

\medskip
If we require that $d_k \le d/2^k$, that is, we find at least half of
the roots at each step, then we can relax the restriction on $\|r_k\|$
somewhat.   If we have $\|r_k\| \le \frac{1}{4\cdot4^d}$, then we can
use Rouch\'e's theorem as before to show that the roots of $p_k f_k$ are in
the disk of radius 
$C_k = \frac{1}{2} + \frac{1}{2}\sum_{n=1}^{k}2^{-2^k} < 3/4$.   Note that
for $|z| = C_k$, we have 
$$ |f_{k-1}(z)| \ge (C_k - C_{k-1})^{d_{k-1}}
\ge \(\frac{1}{2\cdot 2^{2^k}}\)^{d/2^{k-1}} \ge 1/4^d .$$
 Applying
almost the same argument as before, we have
$$ |f_{k-1}(z) - p_k(z) f_k(z)| 
   \le \|r_k\| \sum_{j=0}^{\deg r_k}|z|^j
   \le \frac{1}{4^d}
   \le |f_{k-1}(z)|. $$
Since $\|r_k\| < \mu (4/7)^{d+3} < \frac{1}{4\cdot4^d}$, we are done.
\qquad\end{proof}

\medskip
\subsection{Controlling the error} \Label{err-bound}
In this section, we compute the size of $\tau$ that we can use to ensure we
have an $\epsilon$-factorization of $\psi$. The following very simple
proposition shows that if the norms of two polynomials are close, so are the
norms of the rescaled versions. This gives us the relationship between
$\epsilon$ and the number $\mu$ used in \lemref{repeated-deflation}.

\begin{proposition}
Suppose that $\phi$ is a monic polynomial with all its roots in $\D_{R}$, so
that the roots of $f(z) = (2R)^{-d} \phi(2Rz)$ are in $\D_{1/2}$.  Then if
$\|f - p\| \le \epsilon$, we have
	$$ \left\| \phi(z) - (2R)^d p\(\frac{z}{2R}\) \right\| \le
        (2R)^d\epsilon. $$
\end{proposition}

\begin{proof}
Let $f(z) = \sum a_j z^j$ and $p(z) = \sum b_j z^j$.  Then
$\phi(z)=\sum(2R)^{d-j}a_j z^j$ and $(2R)^d p\(\frac{z}{2R}\) =
 \sum(2R)^{d-j}b_j z^j$.  Since $\max | \sum(2R)^{d-j}(a_j -b_j) | \le
 (2R)^d |a_j - b_j|$, we have our claim.
\qquad\end{proof}

\bigskip Our input polynomial $\phi$ is in $\Pd(1)$ and $f$ is be the
rescaled polynomial as in the previous proposition, so $R=2$.  Then an
$\epsilon$-factorization of $\phi$ corresponds to an
$\epsilon/4^d$-factorization of $f$, so we take $\mu=\epsilon/4^d$.

In order to properly approximate the roots of $f$, we need to ensure that the 
remainder $r_k$ (as in \lemref{repeated-deflation}) at the $k^{th}$ step
satisfies 
$\|r_k\| < 2\tau$, where $\tau = (4/7)^{d+3}\mu/2 = 32\epsilon/7^{d+3}$.

At each stage, we translate $f_k$ by $\tau$, and ensure that the error
introduced by the deflation of the translated polynomial is no more than
$\tau$.  By \lemref{single-deflation}, we need the root distance between
translated polynomial and the deflated polynomial to satisfy
$$ \delta \le
	\frac{\tau}{8d}\(\frac{4}{7}\)^d.$$
Then we will have
$$ \|r_k\| = \|f_{k-1} - p_k f_k\| \le 
	\| (f_{k-1} - \vec\tau) - p_k f_k \| + 
        \| (f_{k-1} - \vec\tau) - f_{k-1} \|
        \le \tau + \tau. $$ 

\medskip
In order to achieve the root distance less than $\delta$, we apply Newton's
method to the approximate zeros found by the routine {\tt iterate-PLM} (see
\lemref{PLM-its}).  Since each point $z$ is an approximate zero to the root
$\xi$ of the 
translated polynomial, we have by \defref{approx-zero}
$$ \delta \le 8\(\frac{1}{2}\)^{2^n}|z-\xi|.$$
Thus, iterating Newton's method $\log_2 \log_2 (8/\delta)$ times, as is done
in {\tt polish-roots},  will give the desired result.

\medskip
\subsection{Arithmetic complexity} \Label{complexity}
In this section we count the number of arithmetic operations involved in
using the algorithm to obtain an $\epsilon$-factorization of a polynomial in
$\Pd(1)$. 

\medskip\noindent
In the main routine, we rescale the polynomial and then invoke 
{\tt get-\dh{}half-\dh{}roots-\dh{}and-deflate} at most $\log_2 d$ times,
since at least half of the roots are found in each call.  

\begin{list}
 {$\bullet$}
 {\setlength{\itemsep}{.25\baselineskip}
  \setlength{\topsep}{.5\baselineskip}
  \setlength{\leftmargin}{1em}
  \setlength{\itemindent}{0pt}
 }
  \item The cost of rescaling is $2d$ multiplications.
  \item Each call to {\tt get-half-roots-and-deflate} calls 
        {\tt choose-\dh{}4d-\dh{}good-\dh{}initial-\dh{}points} once, and
        makes at most 4 calls to each of {\tt iterate-\dh{}PLM}, 
        {\tt select-\dh{}approx-\dh{}zeros}, {\tt polish-\dh{}roots},  
        and {\tt weed-\dh{}out-\dh{}duplicates}, and one call to {\tt deflate}.
        As before, we denote the degree of the input polynomial by $d$ and
        the degree of the $k^{th}$ intermediate polynomial by $d_k$.

  \medskip
\begin{list}
 {$\circ$}
 {\setlength{\itemsep}{.25\baselineskip}
  \setlength{\topsep}{.25\baselineskip}
  \setlength{\leftmargin}{1em}
  \setlength{\itemindent}{0pt}
 }
    \item The routine {\tt choose-4d-good-initial-points} involves
          evaluation of $f_k$ at $676d_k$ points, and $676d_k$ comparisons.
          Since we can evaluate $f_k$ at $m$ points with a cost of $\Order{m
          (\log d_k)^2}$ operations (see \cite{BorodinMunro}), this gives a
          total of $\Order{d_k (\log d_k)^2}$ operations.
    \item For {\tt iterate-PLM}, each iteration evaluates $z_n -
          \frac{f(z_n) - w_{n+1}}{f'(z_n)}$ at $d_k$ points, which costs 
          $\Order{d_k (\log d_k)^2}$.  This is done $N$ times, where 
          $N < 27 \log(\frac{|w_0|}{\tau})$.  Note that since the roots of $f$
          are in $\D_{3/4}$ and $|z_0|=3/2$, we have $|w_0| = |f_k(z_0)|
          \le (9/4)^{d_k}$.  Since $\tau = 32\epsilon/7^{d+3}$, we have 
          $N = C_1(d+|\log\epsilon|)$ for some constant $C_1$.  This
          gives a total of  
          $$\Order{(d + |\log\epsilon|) (d_k(\log d_k)^2)}$$
          operations.
    \item For {\tt select-approx-zeros}, we need to evaluate $\alpha$ at
          $d_k$ points, which requires evaluation of all the derivatives of
          $f$.  This requires $\Order{d_k^2 (\log d_k)^2}$ operations.
    \item The routine {\tt polish-roots} performs $M$ iterations of Newton's
          method at $d_k$ points.   Since $M = C_2(\log d +
          \log|\log\epsilon|)$ for some constant $C_2$, we have a total
          operation count of 
          $$ \Order{(d_k(\log d_k)^2)(\log d + \log |\log\epsilon|)}.$$
    \item {\tt Weed-out-duplicates} requires a sort of at most $d_k$ points,
          which costs $\Order{d_k \log d_k}$, and evaluation of $f$ and $f'$
          at $d_k$ points.  Thus the total cost of this routine is
          $\Order{d_k (\log d_k)^2}$. 
    \item Lastly, {\tt deflate} costs $\Order{d_k(\log d_k)^2}$ operations.
  \end{list}
\end{list}

\medskip\noindent
Thus, the overall cost of each call to {\tt get-half-roots-and-deflate}
is dominated by that of {\tt iterate-PLM}, and is at most 
$\Order{d_k(\log d_k)^2(d + |\log\epsilon|)}$ operations.
Since $d_{k+1} \le d_k/2$, we have a total cost of at most
	$$\Order{d^2(\log d)^2 + d(\log d)^2|\log\epsilon|}$$
operations to obtain an $\epsilon$-factorization of the polynomial.

\medskip
\subsection{Summary and proof of main theorem}
\Label{proof-of-main-thm}
At this point, we have actually already proven \thmref{main-theorem}, but we
would like to tie together the various steps involved.  Just to refresh your
memory, this theorem says, in essence, that our algorithm always produces an
$\epsilon$-factorization with the stated complexity.

Recall that the algorithm performs the approximate factorization in stages
using the routine {\tt get-half-roots-and-deflate};
at the $k^{th}$ step, we produce a function $f_k$ and sets of approximations
$\Lambda_j$ so that
  $$f(z) \approx f_k(z) \prod_{\lambda_i \in \Lambda} (z - \lambda_i),$$
where $\Lambda = \Union_{j\le k}\Lambda_j$.  In order to prove
\thmref{main-theorem}, we need to show that the number of approximations found
at the $k^{th}$ stage ($\#\Lambda_k$) is at least $(\deg f_{k-1})/2$, and
that $\|f(z) - f_k(z) \prod_{\lambda_i \in \Lambda} (z - \lambda_i)\| \le
\epsilon/4^d$. Given this, the complexity calculation in the previous
section will apply, and we shall have the result.

\smallskip
First, note that as a consequence of \lemref{good-sector-count}, there will
always be a quarter-plane in the target space on which there at least 
$d_{k}/2$ branches of $f_{k}^{-1}$ are defined.  This means that if we
start with $d_k$ points $z_j$ which are well-spaced in the source space (so
that each sheet in the target space is represented),
then for at least half of them there will be a branch of the the inverse
$f_{z_j}^{-1}$ which is defined on the entire quadrant, and $f_{z_i}^{-1}
\ne f_{z_j}^{-1}$.  For these points $z_j$, if  we ensure that $f(z_j)$ is
close enough to the center line of the quadrant, at least half of them will
satisfy the hypothesis of \thmref{PLM-follows-ray} and so the routine {\tt
iterate-PLM} will produce approximations to each of the corresponding
$d_k/2$ roots of $f_k$ (with multiplicity).  Such initial points $z_j$ will
be produced by {\tt choose-4d-good-initial-points}, as was shown in
\secref{initial-points}. As a consequence of \lemref{PLM-its}, the
good approximations   {\tt iterate-PLM} are approximate zeros of 
$\psi = f_k - \vec{\tau}$, with $\alpha_\psi < 1/8$.   As was discussed in
\secref{weeding}, application of the routines {\tt select-approx-zeros} and 
{\tt weed-out-duplicates} will select exactly one representative for each
approximated root of $\psi$, giving at least $d_k/2$ such approximations
$\lambda_j$.  This selection is necessary since some of the initial $z_j$
which did not satisfy the hypothesis of \thmref{PLM-follows-ray} may still
have converged. 

As was shown in \secref{err-bound}, the approximations produced yield an
$\epsilon/4^d$-fact\-or\-ization, since the each the $\lambda_i$ are made
sufficiently close to the roots of $\psi$ by the routine {\tt polish-roots},
and $\|\psi - f_k\|$ is sufficiently small.  This completes the proof of
\thmref{main-theorem}.

%
\nocite{Manning}
\bibliography{scott}
\bibliographystyle{siam}

\end{document}